\documentclass[10pt,reqno]{amsart}
\usepackage{amsmath,amssymb,amsthm,mathtools,mathrsfs}
\usepackage{enumitem,needspace}
\usepackage{graphicx}
\usepackage[margin=1in]{geometry}
\calclayout
\usepackage[section]{placeins}
\usepackage{flafter}
\usepackage{cite}
\usepackage[hidelinks]{hyperref}
\usepackage{bookmark}
\numberwithin{equation}{section}

\makeatletter
\renewcommand{\subsection}{\@startsection{subsection}{2}%
  \z@{.5\linespacing\@plus.7\linespacing}{.25\linespacing}%
  {\normalfont\bfseries}}
\makeatother

\theoremstyle{plain}
\newtheorem{theorem}{Theorem}[section]
\newtheorem{proposition}[theorem]{Proposition}
\newtheorem{lemma}[theorem]{Lemma}
\newtheorem{corollary}[theorem]{Corollary}
\theoremstyle{definition}
\newtheorem{definition}{Definition}[section]
\newtheorem{example}[theorem]{Example}
\theoremstyle{remark}
\newtheorem{remark}[theorem]{Remark}
\theoremstyle{plain}

\newcommand{\R}{\mathbb R}

\newcommand{\Ts}{T^{\ast}}
\newcommand{\Cc}{C_c(T)}
\newcommand{\Wa}{W_a^{1,p}(T,m)}
\newcommand{\Waz}{W_{a,0}^{1,p}(T,m)}
\newcommand{\Qa}{Q_{p,a}}
\DeclareMathOperator{\Ran}{Ran}
\newcommand{\ind}{\mathbf 1}
\newcommand{\av}[2]{\langle #1\rangle_{#2}}
\newcommand{\norm}[1]{\lVert #1\rVert}

\DeclareMathOperator{\supp}{supp}
\DeclareMathOperator{\sgn}{sgn}
\DeclareMathOperator{\dist}{dist}
\setlist[enumerate]{label=\textup{(\roman*)},leftmargin=2.4em,itemsep=2pt,topsep=4pt}
\allowdisplaybreaks[1]
\hypersetup{pdftitle={Finite-Support Approximation in Variable-Coefficient Sobolev Spaces on Nonhomogeneous Flow-Measure Trees},pdfauthor={Zeyu Cheng; Yichun Zhao; Jiang Zhou},pdfsubject={Finite-support approximation in variable-coefficient Sobolev spaces on flow-measure trees},pdfkeywords={nonhomogeneous tree, flow measure, variable-coefficient Sobolev space, finite-support approximation}}
\title[Finite-support approximation on flow-measure trees]{\bfseries Finite-Support Approximation in Variable-Coefficient Sobolev Spaces on Nonhomogeneous Flow-Measure Trees}
\author[Z. Cheng]{Zeyu Cheng}
\address{Zeyu Cheng: School of Mathematics and System Sciences, Xinjiang University\\
Urumqi 830046\\
People's Republic of China}
\email{1025729920@qq.com}
\author[Y. Zhao]{Yichun Zhao}
\address{Yichun Zhao: School of Mathematics and Statistics, Yili Normal University\\
Yining 835000\\
People's Republic of China}
\email{zhaoyichun@stu.xju.edu.cn}
\author[J. Zhou]{Jiang Zhou\textsuperscript{*}}
\address{Jiang Zhou: School of Mathematics and System Sciences, Xinjiang University\\
Urumqi 830046\\
People's Republic of China}
\email{zhoujiang@xju.edu.cn}
\thanks{\textsuperscript{*}Corresponding author.}
\date{}

\begin{document}
\begin{abstract}
Let $T$ be a nonhomogeneous rooted tree equipped with a positive flow measure $m$. For $1<p<\infty$, consider the variable-coefficient Sobolev space on the nonhomogeneous flow-measure tree~$T$
\[
W_a^{1,p}(T,m)
=
\left\{
f\in L^p(T,m):
\sum_{x\ne \rho}a(x)m(x)\lvert f(x)-f(\pi(x))\rvert^p<\infty
\right\},
\]
where $a$ is a positive variable coefficient, $\rho$ denotes the root and $\pi(x)$ denotes the parent of $x$. Using an averaging method over descendant levels, we prove the density of finitely supported functions in $W_a^{1,p}(T,m)$ for coefficients satisfying a weighted moment condition. As an application, the relationship between solutions on finite domains and global solutions is clarified for $p=2$.
\end{abstract}
\maketitle

\section{Introduction}\label{sec:introduction}

The purpose of this paper is to establish variable-coefficient Sobolev spaces on nonhomogeneous flow-measure trees and to investigate conditions ensuring the density of finitely supported functions. We first review the relevant research background.

\subsection{Background and motivation}
The study of harmonic analysis on trees was initiated by Cartier in his investigation of harmonic functions on trees \cite{Cartier1972} and was further developed by Fig\`a-Talamanca \cite{FigaTalamancaNebbia1991}. Later, Hebisch introduced the canonical flow measure on homogeneous trees while studying singular integrals and spectral multipliers \cite{HebischSteger2003}. Building on this framework, Levi established $L^p$ estimates for the Riesz transform on homogeneous flow-measure trees \cite{LMSTV2023}. However, branching networks arising in applications often exhibit nonuniform branching structures and distributions that are not fully captured by homogeneous tree models. To study such structures, analysis has been extended to nonhomogeneous trees equipped with flow measures. Several aspects of harmonic analysis on nonhomogeneous flow-measure trees were developed by Levi \cite{LSTV2023}.

Classical Sobolev spaces were introduced in \cite{Sobolev1938} during the study of boundary value problems for elliptic equations. With the development of function space theory, many results have been established for classical Sobolev spaces \cite{MeyersSerrin1964,Moser1971,Aubin1976,BrezisNirenberg1983,Adams1988}. As the underlying spaces considered in Sobolev theory extended to tree structures, Evans related Sobolev embedding problems to weighted function spaces on trees \cite{EvansHarris1993}. Subsequently, Chung investigated Sobolev inequalities on graphs in 1995 \cite{ChungYau1995}. In 2001, Naimark studied the structure of Sobolev spaces on regular metric trees \cite{NaimarkSolomyak2001}. In 2003, Solomyak investigated approximation problems in Sobolev spaces on metric graphs \cite{Solomyak2003}. In 2005, Ostrovskii studied Sobolev spaces defined on graph vertices and their Banach-space properties \cite{Ostrovskii2005}. In 2025, Shao investigated completeness, separability, and reflexivity of Sobolev spaces on locally finite graphs \cite{ShaoYangZhao2025}. For related operator studies on finite graphs, see \cite{SoriaTradacete2016,WangZhao2026,XiangZhao2026}.

In Sobolev theory on trees, the density of finitely supported functions is an important question. In 2009, Maury investigated finite-energy function spaces and finite-support closures on infinite resistive trees \cite{MaurySalortVannier2009}. In 2018, Nicaise studied the density of compactly supported functions in first-order Sobolev spaces on infinite weighted $p$-adic trees \cite{NicaiseSemin2018}. In 2022, Koskela investigated the density of compactly supported functions in Sobolev spaces on regular metric trees \cite{KoskelaNguyenWang2022}. In 2025, Murmann studied finite-support approximation in Sobolev spaces on weighted discrete graphs \cite{MurmannSchmidt2025}. In the same year, Tian showed that finitely supported functions need not be dense in the Sobolev norm on locally finite graphs \cite{TianZhao2025}.

The present work is motivated by the following three considerations.

\begin{itemize}[leftmargin=1.5em,itemsep=5pt,topsep=3pt]
\item As we all know, the infinite resistive-tree model \cite{MaurySalortVannier2009} represents a tree-shaped branching network that generally has a nonuniform branching structure. The difficulty lies in the fact that branching numbers and measure distributions may vary across vertices. Therefore, to describe both complex branching structures and flow-mass distributions, it is necessary to study nonhomogeneous rooted trees equipped with positive flow measures. This provides the primary motivation for our study.

\item Sobolev spaces were introduced in \cite{Sobolev1938}. Their important role in harmonic analysis and partial differential equations owes much to their ability to control both functions and their weak derivatives, thereby providing a functional framework for weak solutions, energy estimates, and variational methods. Rich results have also been obtained for Sobolev spaces on different underlying spaces \cite{Hajlasz1996,NaimarkSolomyak2001,Ostrovskii2005}. Therefore, our second motivation is to develop related Sobolev space theory on nonhomogeneous rooted trees equipped with positive flow measures.

\item We observe that different edges in a branching network may have different response coefficients, and the flow measure alone may not adequately characterize these differences in edge responses. This is because the flow measure only affects the distribution of reference mass across branches and cannot characterize the response properties of different edges. Therefore, to account for such nonuniform edge responses, we need to introduce a positive variable coefficient into the difference energy. This provides the third motivation for our study.
\end{itemize}

\subsection{Basic assumptions}\label{sec:prelim}
Let $T$ be a connected locally finite tree rooted at a vertex $\rho$. We also use $T$ for its vertex set and write
\[
 \Ts=T\setminus\{\rho\}.
\]
For $x\in\Ts$, let $\pi(x)$ be its parent; for $z\in T$, write
\[
 s(z)=\{x\in\Ts:\pi(x)=z\}.
\]
The parent edge $(\pi(x),x)$ is identified with its nonroot endpoint $x$. Let $|x|$ denote the distance from $x$ to $\rho$, and set
\[
 S_n=\{x\in T:|x|=n\},
 \qquad
 B_N=\bigcup_{n=0}^N S_n.
\]
For $z\in T$, let $T_z$ be the descendant subtree rooted at $z$, and put
\[
 C(z,n)=T_z\cap S_n,\qquad n\ge|z|.
\]
Fix
\[
 1<p<\infty,\qquad p'=\frac{p}{p-1}.
\]
All functions are real-valued. For a set $A$, let $\ind_A$ denote its indicator function. Write $\ell^\infty(T)$ for the bounded functions on $T$, with norm $\norm f_\infty=\sup_{x\in T}|f(x)|$. Let $C_c(T)$ denote the finitely supported functions, and write
\[
 c_0(T)=\{f:T\to(-\infty,\infty):\{x\in T:|f(x)|>\varepsilon\}\text{ is finite for every }\varepsilon>0\}.
\]

\begin{definition}\label{def:flow}
Let $m:T\to(0,\infty)$. If
\begin{equation}\label{eq:flow}
 m(z)=\sum_{x\in s(z)}m(x)\qquad(z\in T),
\end{equation}
then $m$ is called a positive flow measure on $T$. Throughout, $m(\rho)=1$.
\end{definition}

For $A\subseteq T$, write $m(A)=\sum_{x\in A}m(x)$.

By iterating the flow conservation identity \eqref{eq:flow}, for every $z\in T$ and $n\ge|z|$ we obtain
\begin{equation}\label{eq:cellmass}
 \sum_{x\in C(z,n)}m(x)=m(z),\qquad n\ge|z|.
\end{equation}
In particular, since $m(\rho)=1$,
\begin{equation*}
 m(S_n)=1,\qquad m(B_N)=N+1,\qquad m(T)=\infty.
\end{equation*}

With the flow measure specified, we consider the corresponding $L^p$ space.

\begin{definition}\label{def:lp}
Let $m$ be a positive flow measure on $T$. If $f:T\to(-\infty,\infty)$ satisfies
\[
 \sum_{x\in T}|f(x)|^pm(x)<\infty,
\]
then $f\in L^p(T,m)$, with
\[
 \norm f_p=\left(\sum_{x\in T}|f(x)|^pm(x)\right)^{1/p}.
\]
The space $L^p(\Ts,m)$ is defined analogously.
\end{definition}

To establish the variable-coefficient Sobolev space, we also need to describe how functions vary across adjacent vertices.

\begin{definition}\label{def:difference}
For $f:T\to(-\infty,\infty)$, define
\begin{equation*}
 Df(x)=f(x)-f(\pi(x)),\qquad x\in\Ts.
\end{equation*}
\end{definition}

We now introduce the definition of the variable-coefficient Sobolev space.

\begin{definition}\label{def:space}
Let $m$ be a positive flow measure, let $1<p<\infty$, and let
\begin{equation*}
 a:\Ts\to(0,\infty)
\end{equation*}
be finite-valued. For $f:T\to(-\infty,\infty)$, define
\begin{equation*}
 Q_{p,a}(f)=\sum_{x\in\Ts}a(x)m(x)|Df(x)|^p.
\end{equation*}
The corresponding variable-coefficient Sobolev space is defined by
\begin{equation*}
 W_a^{1,p}(T,m)=\{f\in L^p(T,m):Q_{p,a}(f)<\infty\},
\end{equation*}
equipped with
\begin{equation*}
 \norm f_{W_a^{1,p}}=(\norm f_p^p+Q_{p,a}(f))^{1/p}.
\end{equation*}
\end{definition}

To study the density of finitely supported functions in this space, we consider their closure in the Sobolev norm.

\begin{definition}\label{def:minimal}
Let $W_a^{1,p}(T,m)$ be as in Definition~\ref{def:space}. Define
\begin{equation*}
 W_{a,0}^{1,p}(T,m)=\overline{C_c(T)}^{\,W_a^{1,p}(T,m)}.
\end{equation*}
Here the closure is taken with respect to the $W_a^{1,p}(T,m)$ norm.
\end{definition}

We next introduce descendant-level averages and a weighted moment constant as tools for studying finite-support approximation.

For $z\in T$, $n\ge |z|$, and any real-valued function $h$ on $C(z,n)$, write its descendant-level average as
\[
 \av{h}{z,n}=\frac{1}{m(z)}\sum_{x\in C(z,n)}m(x)h(x).
\]

Write the weighted moment constant as
\begin{equation*}
 \mathfrak C_p(a)=\sup_{\substack{z\in T\\n>|z|}}
 \av{a}{z,n}\av{a^{-1/(p-1)}}{z,n}^{p-1}.
\end{equation*}

\subsection{Main theorems and contributions}\label{sec:intro-main}
We now present the main theorems of this paper. The first theorem establishes the basic structural properties of the variable-coefficient Sobolev space.

\begin{theorem}\label{thm:structure}
The space $W_a^{1,p}(T,m)$ is a separable, uniformly convex Banach space and hence reflexive. When $p=2$, it is a Hilbert space with inner product
\[
 \begin{aligned}
 (f,g)_{W_a^{1,2}}
 &=\sum_{x\in T}f(x)g(x)m(x)\\
 &\quad+\sum_{x\in T^*}a(x)m(x)Df(x)Dg(x).
 \end{aligned}
\]
\end{theorem}

\Needspace{9\baselineskip}
We next investigate the density of finitely supported functions in this space.

\begin{theorem}\label{thm:density}
If $\mathfrak C_p(a)<\infty$, then for every $f\in W_a^{1,p}(T,m)$, there exists a sequence $\{f_j\}_{j\ge1}\subset C_c(T)$ such that
\begin{equation*}
 \|f_j-f\|_{W_a^{1,p}(T,m)}\longrightarrow 0
 \qquad (j\to\infty).
\end{equation*}
Moreover, $C_c(T)$ is a dense subspace of $W_a^{1,p}(T,m)$.
\end{theorem}

\Needspace{11\baselineskip}
The density of finitely supported functions is also closely related to the convergence of finite-domain solutions to the global variational solution. We now describe this relationship for $p=2$.

\begin{theorem}\label{thm:finite-domain-convergence}
Let $p=2$ and use the notation in the following remark. For any positive, finite-valued coefficient $a$, the finite-domain solutions $u_N$ satisfy
\[
 \|u_N-u^{\max}\|_H\longrightarrow0\qquad(N\to\infty)
\]
for every $h\in L^2(T,m)$ if and only if $C_c(T)$ is dense in $H$.
\end{theorem}

\Needspace{7\baselineskip}
\noindent\textbf{Remark.}\ By Theorem~\ref{thm:structure}, when $p=2$, the variable-coefficient Sobolev space
\[
 H:=W_a^{1,2}(T,m)
\]
is a Hilbert space. Denote its inner product by
\[
 B_a(u,v):=(u,v)_H.
\]
For each nonnegative integer $N$, let
\[
 V_N=\{v\in H:\supp v\subseteq B_N\}
\]
be the corresponding finite-dimensional subspace of $H$.
For each $h\in L^2(T,m)$, let $u^{\max}\in H$ and $u_N\in V_N$ denote the unique solutions of the following variational problems, respectively:
\[
 \begin{aligned}
 B_a(u^{\max},v)&=\langle h,v\rangle_{L^2(T,m)},&&v\in H,\\
 B_a(u_N,v)&=\langle h,v\rangle_{L^2(T,m)},&&v\in V_N.
 \end{aligned}
\]

\Needspace{8\baselineskip}
The main contributions are as follows.

\begin{itemize}[leftmargin=1.5em,itemsep=5pt,topsep=3pt]
\item We establish and study the variable-coefficient Sobolev space $W_a^{1,p}(T,m)$ on nonhomogeneous flow-measure trees. When $a\equiv1$, this space has the same function domain as $L^p(T,m)$, and the norms of the two spaces are equivalent; see Lemma~\ref{lem:D-bound} and Proposition~\ref{prop:lp-equality}. Our approach differs from classical Sobolev theory on Euclidean domains. The key obstacle is that classical differentiation methods cannot be directly applied to functions defined on the vertices of graphs and trees. We therefore use differences between adjacent vertices to describe variations of functions.

\item We construct descendant-level averaging operators and establish the density of finitely supported functions under a weighted moment condition. This approach differs from that in \cite{KoskelaNguyenWang2022}. The main difficulty is that branching numbers and flow distributions may vary across vertices on nonhomogeneous flow-measure trees, while variable coefficients further complicate the control of difference energies.

\item Our Theorem~\ref{thm:density} shows that, under the weighted moment condition, variational identities that are continuous in the Sobolev norm can be verified on finitely supported test functions and then extended to the entire Sobolev space. This allows variational problems on infinite trees to be studied through finitely supported tests, even when the coefficients vary across branches.

\item Our Theorem~\ref{thm:finite-domain-convergence} investigates the relationship between finite-domain solutions and global variational solutions in a diffusion model on branching media. For $p=2$, finite-domain solutions converge to the global variational solution in the $W_a^{1,2}(T,m)$ norm for every $h\in L^2(T,m)$ if and only if finitely supported functions are dense in $W_a^{1,2}(T,m)$.
\end{itemize}

\section{Proofs of the main results}\label{sec:density}

\subsection{Basic structure of the space}
\begin{lemma}\label{lem:D-bound}
For every $f\in L^p(T,m)$,
\begin{equation}\label{eq:D-bound}
 \norm{Df}_{L^p(\Ts,m)}\le2\norm f_p.
\end{equation}
\end{lemma}
\begin{proof}
By the flow identity \eqref{eq:flow},
\begin{equation*}
 \begin{aligned}
 \sum_{x\in\Ts}|f(\pi(x))|^pm(x)
 &=\sum_{z\in T}|f(z)|^p\sum_{x\in s(z)}m(x)\\
 &=\norm f_p^p.
 \end{aligned}
\end{equation*}
Since $Df(x)=f(x)-f(\pi(x))$, the triangle inequality and
$\norm{f|_{\Ts}}_{L^p(\Ts,m)}\le\norm f_p$ yield \eqref{eq:D-bound}.
\end{proof}

\noindent\textbf{Proof of Theorem~\ref{thm:structure}.} We now prove Theorem~\ref{thm:structure}.
\begin{proof}
For $f\in W_a^{1,p}(T,m)$, let $G_af=a^{1/p}Df$. We first show that $G_a$, viewed as an operator from $L^p(T,m)$ to $L^p(\Ts,m)$ with domain $W_a^{1,p}(T,m)$, is closed.

Suppose $(f_j)\subset W_a^{1,p}(T,m)$ satisfies
\[
 f_j\to f\quad\text{in }L^p(T,m),\qquad
 G_af_j\to g\quad\text{in }L^p(\Ts,m).
\]
Since $m(x)>0$ at every vertex, these convergences imply pointwise convergence. Hence, for every $x\in\Ts$,
\[
 \begin{aligned}
 g(x)&=\lim_{j\to\infty}a(x)^{1/p}
 \bigl(f_j(x)-f_j(\pi(x))\bigr)\\
 &=a(x)^{1/p}\bigl(f(x)-f(\pi(x))\bigr)=G_af(x).
 \end{aligned}
\]
Consequently, $f\in W_a^{1,p}(T,m)$ and $G_af=g$. Thus $G_a$ is closed.

Consider the linear map
\begin{equation*}
 J:W_a^{1,p}(T,m)\longrightarrow L^p(T,m)\oplus_p L^p(\Ts,m),
 \qquad Jf=(f,G_af),
\end{equation*}
where $\oplus_p$ denotes the direct sum endowed with the $\ell^p$ norm. By the definition of the Sobolev norm,
\[
 \norm{Jf}^p=\norm f_p^p+\norm{G_af}_{L^p(\Ts,m)}^p
 =\norm f_{W_a^{1,p}}^p.
\]
Thus $J$ is a linear isometry. Since $G_a$ is closed, its graph $J(W_a^{1,p}(T,m))$ is a closed linear subspace of the $L^p$ direct sum.

For $1<p<\infty$, this direct sum is a separable, uniformly convex Banach space. Its closed subspace $J(W_a^{1,p}(T,m))$ therefore has the same properties. Hence $W_a^{1,p}(T,m)$ is a separable, uniformly convex Banach space and thus reflexive. When $p=2$, the direct sum is a Hilbert space, and the induced inner product is precisely the one in Theorem~\ref{thm:structure}.
\end{proof}

\begin{proposition}\label{prop:lp-equality}
The following conditions are equivalent:
\begin{enumerate}
\item $Q_{p,a}(f)<\infty$ for every $f\in L^p(T,m)$;
\item $G_a$ is defined on all of $L^p(T,m)$ and is bounded as an operator into $L^p(\Ts,m)$;
\item $\sup_{x\in\Ts}a(x)<\infty$.
\end{enumerate}
Under these conditions, every $f\in L^p(T,m)$ satisfies
\begin{equation}\label{eq:equiv-lp}
 \norm f_p\le\norm f_{W_a^{1,p}}
 \le(1+2^p\norm a_\infty)^{1/p}\norm f_p,
\end{equation}
where $\norm a_\infty=\sup_{x\in\Ts}a(x)$.
\end{proposition}
\begin{proof}
Suppose that $a$ is bounded. By Lemma~\ref{lem:D-bound},
\[
 \begin{aligned}
 \Qa(f)&\le\norm a_\infty\norm{Df}_{L^p(\Ts,m)}^p\\
 &\le2^p\norm a_\infty\norm f_p^p.
 \end{aligned}
\]
This proves (iii)$\Rightarrow$(ii)$\Rightarrow$(i), together with \eqref{eq:equiv-lp}.

Conversely, suppose that (i) holds. As shown in the proof of Theorem~\ref{thm:structure}, $G_a$ is a closed operator. Since its domain is now all of $L^p(T,m)$, the closed graph theorem yields a constant $C>0$ such that
\[
 \Qa(f)\le C^p\norm f_p^p.
\]
For each $x\in\Ts$, taking $f=\ind_{\{x\}}$ gives
\begin{equation*}
 \begin{aligned}
 \Qa(\ind_{\{x\}})&=a(x)m(x)+\sum_{y\in s(x)}a(y)m(y),\\
 \norm{\ind_{\{x\}}}_p^p&=m(x).
 \end{aligned}
\end{equation*}
Hence $a(x)m(x)\le C^pm(x)$. Since $m(x)>0$, we obtain $a(x)\le C^p$ for every $x\in\Ts$, proving (iii).
\end{proof}

Proposition~\ref{prop:lp-equality} shows that if $a$ is unbounded, there exists $f\in L^p(T,m)$ with $Q_{p,a}(f)=\infty$. On the other hand, every finitely supported function has finite weighted difference energy. Thus the density of $C_c(T)$ in the ambient space $L^p(T,m)$ does not guarantee density in the Sobolev norm.

\subsection{Descendant-level averaging and finite-support approximation}
By \eqref{eq:cellmass}, the descendant-level average $\av{h}{z,n}$ introduced in Subsection~\ref{sec:prelim} is taken with respect to a probability measure.

\begin{lemma}\label{lem:weighted-average}
Let $F$ be a finite set equipped with a probability measure $\mu$ assigning positive mass to every point, and let $w:F\to(0,\infty)$. The averaging operator
\[
 \mathsf Eh=\left(\int_Fh\,d\mu\right)\ind_F
\]
on $L^p(F,w\,d\mu)$ satisfies
\begin{equation}\label{eq:averaging-norm}
 \norm{\mathsf E}^p
 =\left(\int_Fw\,d\mu\right)
  \left(\int_Fw^{-1/(p-1)}\,d\mu\right)^{p-1}.
\end{equation}
\end{lemma}
\begin{proof}
By H\"older's inequality,
\[
 \left|\int_Fh\,d\mu\right|^p
 \le\left(\int_Fw|h|^p\,d\mu\right)
     \left(\int_Fw^{-1/(p-1)}\,d\mu\right)^{p-1}.
\]
Multiplying by $\int_Fw\,d\mu$ gives the upper bound. Equality is attained for $h=w^{-1/(p-1)}$, which proves \eqref{eq:averaging-norm}.
\end{proof}

For $n>|z|$, define
\begin{equation*}
 K_p(a;z,n)=\av{a}{z,n}\av{a^{-1/(p-1)}}{z,n}^{p-1},
 \qquad
 \mathfrak C_p(a)=\sup_{z\in T,\ n>|z|}K_p(a;z,n).
\end{equation*}
By H\"older's inequality, each $K_p(a;z,n)$ is finite and satisfies $K_p(a;z,n)\ge1$, with equality if and only if $a$ is constant on $C(z,n)$. By Lemma~\ref{lem:weighted-average}, $K_p(a;z,n)$ is exactly the $p$th power of the operator norm of the $m$-averaging operator on $L^p(C(z,n),a\,dm)$.

For each integer $N\ge0$ and $z\in S_N$, set $F_{z,n}=\av{f}{z,n}$ for $n\ge N$. Since $C(z,N)=\{z\}$, we have $F_{z,N}=f(z)$. Define
\begin{equation*}
 P_Nf(x)=
 \begin{cases}
 f(x),&|x|\le N,\\
 F_{z,|x|},&|x|>N,\quad x\in T_z,\ z\in S_N.
 \end{cases}
\end{equation*}
For every $x$ with $|x|>N$, there is a unique ancestor $z\in S_N$, so $P_N$ is well defined. It leaves $f$ unchanged on $B_N$ and replaces its values on each descendant level $C(z,n)$, $n>N$, by their $m$-average. In general, $P_Nf$ need not have finite support.

For functions $h:\Ts\to(-\infty,\infty)$, define the corresponding edge-averaging operator
\begin{equation}\label{eq:edge-average}
 E_Nh(x)=
 \begin{cases}
 h(x),&1\le|x|\le N,\\
 \av{h}{z,|x|},&|x|>N,\quad x\in T_z,\ z\in S_N.
 \end{cases}
\end{equation}
Both operators are defined by finite sums. Since $E_N$ acts separately on each level and every level is finite, it maps $C_c(\Ts)$ into itself.

\begin{lemma}\label{lem:compatibility}
For every real-valued $f$ and $N\ge0$,
\begin{equation}\label{eq:intertwining}
 D(P_Nf)=E_N(Df).
\end{equation}
\end{lemma}
\begin{proof}
For $z\in S_N$ and $n>N$, flow conservation gives
\begin{align*}
 \sum_{x\in C(z,n)}m(x)f(\pi(x))
 &=\sum_{y\in C(z,n-1)}f(y)\sum_{x\in s(y)}m(x)\\
 &=\sum_{y\in C(z,n-1)}m(y)f(y).
\end{align*}
Hence
\begin{equation*}
 F_{z,n}-F_{z,n-1}=\av{Df}{z,n}.
\end{equation*}
At the interface $n=N+1$, the identity remains valid because $C(z,N)=\{z\}$ and
\[
 F_{z,N+1}-f(z)
 =\frac1{m(z)}\sum_{x\in s(z)}m(x)\bigl(f(x)-f(z)\bigr).
\]
Thus \eqref{eq:intertwining} holds on all edges of depth greater than $N$. On edges of depth at most $N$, both sides equal $Df$.
\end{proof}

To control the weighted difference energy under descendant-level averaging, consider the weighted $L^p$ space
\[
 Y_a=L^p(\Ts,a\,dm),\qquad
 \norm h_{Y_a}^p=\sum_{x\in\Ts}a(x)m(x)|h(x)|^p.
\]
Then $Q_{p,a}(f)=\norm{Df}_{Y_a}^p$.

\begin{proposition}\label{prop:edge-operator}
The condition $\mathfrak C_p(a)<\infty$ is equivalent to uniform boundedness of the family $\{E_N\}_{N\ge0}$ on $C_c(\Ts)$ with respect to the $Y_a$ norm. Under this condition, each $E_N$ extends to a bounded linear operator on $Y_a$, and
\begin{equation}\label{eq:edge-exact-norm}
 \norm{E_N}_{Y_a\to Y_a}^p
 =\sup_{\substack{z\in S_N\\n>N}}K_p(a;z,n),
 \qquad\sup_{N\ge0}\norm{E_N}_{Y_a\to Y_a}^p=\mathfrak C_p(a).
\end{equation}
\end{proposition}
\begin{proof}
Fix $N\ge0$, and put
\[
 M_N=\sup_{\substack{z\in S_N\\n>N}}K_p(a;z,n).
\]
The blocks $C(z,n)$, where $z\in S_N$ and $n>N$, together with the singleton blocks of depth at most $N$, partition $\Ts$. By Lemma~\ref{lem:weighted-average}, for every $h\in C_c(\Ts)$,
\[
 \begin{aligned}
 &\sum_{x\in C(z,n)}a(x)m(x)|E_Nh(x)|^p\\
 &\quad\le K_p(a;z,n)
 \sum_{x\in C(z,n)}a(x)m(x)|h(x)|^p.
 \end{aligned}
\]
On each singleton block, $E_N$ acts as the identity. Since $K_p(a;z,n)\ge1$, summing over all blocks yields
\[
 \norm{E_Nh}_{Y_a}^p\le M_N\norm h_{Y_a}^p.
\]
This estimate also holds for $N=0$, when there are no singleton blocks.

Conversely, fix $z\in S_N$ and $n>N$, and choose
\[
 h(x)=\begin{cases}
 a(x)^{-1/(p-1)},&x\in C(z,n),\\
 0,&\text{otherwise}.
 \end{cases}
\]
By the equality case in Lemma~\ref{lem:weighted-average}, this finitely supported function satisfies
\[
 \frac{\norm{E_Nh}_{Y_a}^p}{\norm h_{Y_a}^p}=K_p(a;z,n).
\]
Consequently,
\[
 \sup_{\substack{h\in C_c(\Ts)\\h\ne0}}
 \frac{\norm{E_Nh}_{Y_a}^p}{\norm h_{Y_a}^p}=M_N.
\]
Every descendant block $C(z,n)$ with $n>|z|$ occurs when $N=|z|$. Hence $\sup_{N\ge0}M_N=\mathfrak C_p(a)$, proving the asserted equivalence for uniform boundedness.

If $\mathfrak C_p(a)<\infty$, the density of $C_c(\Ts)$ in $Y_a$ gives a unique bounded extension of each $E_N$, and \eqref{eq:edge-exact-norm} follows. Since $a(x)m(x)>0$ at every point, convergence in $Y_a$ implies pointwise convergence, so the extensions agree with \eqref{eq:edge-average}.
\end{proof}

The equivalence above concerns uniform boundedness of the edge-averaging family on $Y_a$. For the density of $C_c(T)$ in $W_a^{1,p}(T,m)$, Theorem~\ref{thm:density} gives $\mathfrak C_p(a)<\infty$ as a sufficient condition.

\begin{lemma}\label{lem:PN-bounds}
Suppose $\mathfrak C_p(a)=K<\infty$. Then, for every $f\in W_a^{1,p}(T,m)$ and $N\ge0$,
\[
 P_Nf\in W_a^{1,p}(T,m)\cap c_0(T),\qquad
 \norm{P_Nf}_{W_a^{1,p}}\le K^{1/p}\norm f_{W_a^{1,p}}.
\]
Moreover,
\[
 \norm{P_Nf-f}_{W_a^{1,p}}\longrightarrow0\qquad(N\to\infty).
\]
\end{lemma}
\begin{proof}
We first establish the norm bounds. Fix $N\ge0$. For $z\in S_N$ and $n>N$, Jensen's inequality gives
\[
 m(z)|F_{z,n}|^p\le\sum_{x\in C(z,n)}m(x)|f(x)|^p.
\]
Summing over $z\in S_N$, we obtain
\begin{equation}\label{eq:PN-lp-level}
 \sum_{x\in S_n}m(x)|P_Nf(x)|^p
 \le\sum_{x\in S_n}m(x)|f(x)|^p.
\end{equation}
By Lemmas~\ref{lem:compatibility} and~\ref{lem:weighted-average},
\[
 \begin{aligned}
 &\sum_{x\in C(z,n)}a(x)m(x)|D(P_Nf)(x)|^p\\
 &=m(z)\av{a}{z,n}|\av{Df}{z,n}|^p\\
 &\le K_p(a;z,n)
 \sum_{x\in C(z,n)}a(x)m(x)|Df(x)|^p.
 \end{aligned}
\]
Using $K_p(a;z,n)\le K$ and summing over $z\in S_N$, we obtain
\begin{equation}\label{eq:PN-q-level}
 \sum_{x\in S_n}a(x)m(x)|D(P_Nf)(x)|^p
 \le K\sum_{x\in S_n}a(x)m(x)|Df(x)|^p.
\end{equation}
At depths at most $N$, the values and incoming differences are unchanged. Combining the last two estimates with $K\ge1$ gives
\begin{equation*}
 \norm{P_Nf}_p\le\norm f_p,\qquad
 Q_{p,a}(P_Nf)\le KQ_{p,a}(f),\qquad
 \norm{P_Nf}_{W_a^{1,p}}\le K^{1/p}\norm f_{W_a^{1,p}}.
\end{equation*}
In particular, $P_Nf\in W_a^{1,p}(T,m)$.

Next, we prove strong convergence. Both $f-P_Nf$ and its incoming differences vanish at depths at most $N$. The inequality $|r-s|^p\le2^{p-1}(|r|^p+|s|^p)$, together with \eqref{eq:PN-lp-level} and \eqref{eq:PN-q-level}, gives
\[
 \norm{f-P_Nf}_p^p\le2^p\sum_{|x|>N}m(x)|f(x)|^p
\]
and
\[
 Q_{p,a}(f-P_Nf)\le2^{p-1}(1+K)
 \sum_{|x|>N}a(x)m(x)|Df(x)|^p.
\]
Consequently,
\begin{equation*}
\begin{split}
 \norm{f-P_Nf}_{W_a^{1,p}}^p
 &\le2^p\sum_{|x|>N}m(x)|f(x)|^p\\
 &\quad+2^{p-1}(1+K)
 \sum_{|x|>N}a(x)m(x)|Df(x)|^p
 \longrightarrow0.
\end{split}
\end{equation*}
Both tails tend to zero because $f\in W_a^{1,p}(T,m)$.

Finally, we show that $P_Nf\in c_0(T)$. Fix $N\ge0$ and $z\in S_N$. Jensen's inequality gives
\[
 m(z)\sum_{n=N}^{\infty}|F_{z,n}|^p
 \le\sum_{x\in T_z}m(x)|f(x)|^p<\infty.
\]
Since $m(z)>0$, we have $F_{z,n}\to0$ as $n\to\infty$. Given $\varepsilon>0$, for each $z\in S_N$ choose $R_z\ge N$ so that $|F_{z,n}|\le\varepsilon$ whenever $n>R_z$. Since $S_N$ is finite, $R=\max_{z\in S_N}R_z$ is finite, and
\[
 \{x\in T:|P_Nf(x)|>\varepsilon\}\subseteq B_R.
\]
Thus $P_Nf\in c_0(T)$.
\end{proof}

\begin{remark}\label{prop:nested}
Suppose $\mathfrak C_p(a)<\infty$. For $N\ge0$, let $\mathcal R_N$ be the subspace of $W_a^{1,p}(T,m)$ consisting of functions that are constant on every $C(z,n)$, where $z\in S_N$ and $n>N$.

By the definition of $P_N$, we have $P_Nf\in\mathcal R_N$ and $P_Ng=g$ for every $g\in\mathcal R_N$. Hence
\[
 P_N^2=P_N,\qquad
 \mathcal R_N=\Ran(P_N)=\ker(I-P_N).
\]
By Lemma~\ref{lem:PN-bounds}, $P_N$ is a bounded projection, and $\mathcal R_N$ is a closed subspace of $W_a^{1,p}(T,m)$ contained in $c_0(T)$.

For $M\ge N$, the partition associated with level $M$ refines that of level $N$. Consequently,
\[
 \mathcal R_N\subseteq\mathcal R_M,\qquad P_MP_N=P_N.
\]
On the other hand, for $z\in S_N$ and $k>M$, flow conservation gives
\[
 \begin{aligned}
 \av{P_Mf}{z,k}
 &=\frac1{m(z)}\sum_{w\in C(z,M)}m(w)\av{f}{w,k}\\
 &=\av{f}{z,k}.
 \end{aligned}
\]
For $N<k\le M$, the values are unchanged by $P_M$. Thus $P_NP_M=P_N$. Therefore, for every $M\ge N$,
\begin{equation*}
 P_N^2=P_N,\qquad P_NP_M=P_MP_N=P_N,\qquad
 \mathcal R_N\subseteq\mathcal R_M.
\end{equation*}
Finally, the strong convergence $P_Nf\to f$ from Lemma~\ref{lem:PN-bounds} shows that $\bigcup_{N\ge0}\mathcal R_N$ is dense in $W_a^{1,p}(T,m)$.
\end{remark}

We next use a standard contraction argument to approximate functions in $W_a^{1,p}(T,m)\cap c_0(T)$ by finitely supported functions; compare the related contraction arguments in \cite[Lemma~3.1]{MurmannSchmidt2025} and \cite[Lemma~2.3]{SchillingUemura2012}.

\begin{lemma}\label{lem:threshold}
For $t\in\R$, write $t_+=\max\{t,0\}$. For $\varepsilon>0$, let
\[
 \tau_\varepsilon(t)=\sgn(t)(|t|-\varepsilon)_+,\qquad
 r_\varepsilon(t)=t-\tau_\varepsilon(t)=\max\{-\varepsilon,\min\{t,\varepsilon\}\}.
\]
If $g\in W_a^{1,p}(T,m)\cap c_0(T)$, then $\tau_\varepsilon\circ g\in C_c(T)$ and
\begin{equation}\label{eq:threshold-convergence}
 \norm{g-\tau_\varepsilon\circ g}_{W_a^{1,p}}\longrightarrow0
 \qquad(\varepsilon\downarrow0).
\end{equation}
\end{lemma}
\begin{proof}
Since $g\in c_0(T)$, the set $\{x\in T:|g(x)|>\varepsilon\}$ is finite. As
\[
 \supp(\tau_\varepsilon\circ g)
 \subseteq\{x\in T:|g(x)|>\varepsilon\},
\]
we have $\tau_\varepsilon\circ g\in C_c(T)$.

Set $u_\varepsilon=g-\tau_\varepsilon\circ g=r_\varepsilon\circ g$. Since $r_\varepsilon$ is $1$-Lipschitz and $|r_\varepsilon(t)|\le|t|$, we obtain
\[
 |u_\varepsilon(x)|^p\le|g(x)|^p,\qquad
 |Du_\varepsilon(x)|^p\le|Dg(x)|^p.
\]
Both left-hand sides converge pointwise to zero as $\varepsilon\downarrow0$. Since $g\in W_a^{1,p}(T,m)$, dominated convergence in the vertex sum and the weighted edge sum yields $\norm{u_\varepsilon}_{W_a^{1,p}}\to0$, proving \eqref{eq:threshold-convergence}.
\end{proof}

\noindent\textbf{Proof of Theorem~\ref{thm:density}.} We now prove Theorem~\ref{thm:density}.
\begin{proof}
By Lemma~\ref{lem:PN-bounds}, every $f\in W_a^{1,p}(T,m)$ satisfies $P_Nf\to f$ in $W_a^{1,p}(T,m)$ and $P_Nf\in c_0(T)$. For each fixed $j\ge1$, Lemma~\ref{lem:threshold} allows us to choose $0<\varepsilon_j<1/j$ such that
\[
 \norm{P_jf-\tau_{\varepsilon_j}\circ(P_jf)}_{W_a^{1,p}}<1/j.
\]
Set $f_j=\tau_{\varepsilon_j}\circ(P_jf)$. Then $f_j\in C_c(T)$ and
\begin{equation*}
 \begin{aligned}
 \norm{f-f_j}_{W_a^{1,p}}
 &\le\norm{f-P_jf}_{W_a^{1,p}}+\norm{P_jf-f_j}_{W_a^{1,p}}\\
 &<\norm{f-P_jf}_{W_a^{1,p}}+1/j\longrightarrow0.
 \end{aligned}
\end{equation*}
Therefore, $C_c(T)$ is dense in $W_a^{1,p}(T,m)$.
\end{proof}

The descendant-level averaging construction is related to existing results. On regular trees with uniformly split flow, $P_N$ coincides for $N\ge1$ with the radialization operator in \cite[Definition~6]{CohenEtAl2016}, while a related vertex--edge averaging identity appears in \cite[proof of Theorem~4.8]{Mugnolo2013}. The weighted moment condition here is defined over descendant levels, in contrast to the flow-tree $A_p$ condition based on admissible trapezoids in \cite{OSV2025}.

\subsection{Stability and applications of the moment criterion}\label{subsec:moment-applications}
Theorem~\ref{thm:density} provides a sufficient condition for finite-support density. We first consider its stability under coefficient changes and several coefficient classes covered by this condition.

\begin{proposition}\label{prop:comparable}
Let $b:\Ts\to(0,\infty)$ be another coefficient, and suppose
\[
0<c_-b(x)\le a(x)\le c_+b(x),\qquad x\in\Ts,
\]
where $c_-$ and $c_+$ are fixed lower and upper comparison constants, respectively, independent of $x$, with $0<c_-\le c_+<\infty$.
Then, for every function $f:T\to(-\infty,\infty)$,
\[
c_-Q_{p,b}(f)\le Q_{p,a}(f)\le c_+Q_{p,b}(f).
\]
Consequently, for every $f\in L^p(T,m)$,
\[
Q_{p,a}(f)<\infty
\iff
Q_{p,b}(f)<\infty.
\]
On the common Sobolev domain, the norms satisfy
\[
\begin{aligned}
\min\{1,c_-\}\norm{f}_{W_b^{1,p}}^p
&\le \norm{f}_{W_a^{1,p}}^p\\
&\le \max\{1,c_+\}\norm{f}_{W_b^{1,p}}^p.
\end{aligned}
\]
Thus the two Sobolev norms are equivalent, and $C_c(T)$ is dense in $W_a^{1,p}(T,m)$ if and only if it is dense in $W_b^{1,p}(T,m)$. In particular, changing a positive, finite-valued coefficient on finitely many edges does not affect finite-support density.
\end{proposition}
\begin{proof}
Multiplying the coefficient inequalities by $m(x)|Df(x)|^p$ and summing over $x\in\Ts$ yields the energy estimates.
\Needspace{8\baselineskip}
On the common Sobolev domain,
\[
\begin{aligned}
\norm{f}_{W_a^{1,p}}^p
&\ge \norm{f}_p^p+c_-Q_{p,b}(f)\\
&\ge \min\{1,c_-\}\norm{f}_{W_b^{1,p}}^p,\\
\norm{f}_{W_a^{1,p}}^p
&\le \norm{f}_p^p+c_+Q_{p,b}(f)\\
&\le \max\{1,c_+\}\norm{f}_{W_b^{1,p}}^p.
\end{aligned}
\]
Hence the norms are equivalent, so the finite-support closures coincide and finite-support density is preserved. If the coefficients differ on only finitely many edges, their ratios are positive and finite on those edges and equal to $1$ elsewhere. They therefore satisfy the uniform comparability condition.
\end{proof}

\begin{proposition}\label{prop:radial-invariance}
Let $a:\Ts\to(0,\infty)$ be a positive, finite-valued coefficient, and let $(A_n)_{n\ge1}$ satisfy $0<A_n<\infty$. Define
\[
\widetilde a(x)=A_{|x|}a(x),\qquad x\in\Ts.
\]
Then, for every $z\in T$ and $n>|z|$,
\begin{equation}\label{eq:radial-invariance}
K_p(\widetilde a;z,n)=K_p(a;z,n),\qquad
\mathfrak C_p(\widetilde a)=\mathfrak C_p(a).
\end{equation}
In particular, if $\mathfrak C_p(a)<\infty$, then $C_c(T)$ is dense in both $W_a^{1,p}(T,m)$ and $W_{\widetilde a}^{1,p}(T,m)$.
\end{proposition}
\begin{proof}
Since $A_{|x|}=A_n$ is constant on $C(z,n)$,
\[
\av{\widetilde a}{z,n}=A_n\av{a}{z,n},\qquad
\av{\widetilde a^{-1/(p-1)}}{z,n}
=A_n^{-1/(p-1)}\av{a^{-1/(p-1)}}{z,n}.
\]
The factors involving $A_n$ cancel in $K_p$, giving $K_p(\widetilde a;z,n)=K_p(a;z,n)$. Taking suprema proves \eqref{eq:radial-invariance}. If the common moment constant is finite, Theorem~\ref{thm:density} applied to each coefficient yields finite-support density in both spaces.
\end{proof}

\begin{corollary}\label{cor:radial}
Let
\[
a(x)=A_{|x|},\qquad x\in\Ts,
\]
where $(A_n)_{n\ge1}$ is an arbitrary sequence of positive finite numbers. Then
\[
\mathfrak C_p(a)=1.
\]
Consequently, $C_c(T)$ is dense in $W_a^{1,p}(T,m)$. The sequence $(A_n)_{n\ge1}$ may take arbitrary positive finite values, allowing unrestricted growth, decay, and oscillation.
\end{corollary}
\begin{proof}
Applying Proposition~\ref{prop:radial-invariance} to the constant coefficient $1$ gives $\mathfrak C_p(a)=\mathfrak C_p(1)=1$. The density conclusion follows from Theorem~\ref{thm:density}.
\end{proof}

\begin{corollary}\label{cor:scaling}
Let $a:\Ts\to(0,\infty)$ be a positive, finite-valued coefficient, and let $(A_n)_{n\ge1}$ be an arbitrary sequence of positive finite numbers. Suppose there exist fixed constants $0<c_-\le c_+<\infty$ such that
\[
c_-A_{|x|}\le a(x)\le c_+A_{|x|},\qquad x\in\Ts.
\]
Then
\[
\mathfrak C_p(a)\le\frac{c_+}{c_-}<\infty.
\]
Consequently, $C_c(T)$ is dense in $W_a^{1,p}(T,m)$.
\end{corollary}
\begin{proof}
Fix $z\in T$ and $n>|z|$. Since $A_{|x|}=A_n$ on $C(z,n)$, the coefficient bounds imply
\[
\av{a}{z,n}\le c_+A_n,\qquad
\av{a^{-1/(p-1)}}{z,n}\le(c_-A_n)^{-1/(p-1)}.
\]
It follows that
\[
K_p(a;z,n)\le c_+A_n(c_-A_n)^{-1}
=\frac{c_+}{c_-}.
\]
Taking suprema yields $\mathfrak C_p(a)\le c_+/c_-$. The density conclusion follows from Theorem~\ref{thm:density}.
\end{proof}

\begin{remark}\label{rem:radial-scope}
Proposition~\ref{prop:radial-invariance} and Corollaries~\ref{cor:radial}--\ref{cor:scaling} illustrate how radial invariance of the moment constant yields density results. This invariance, however, does not imply that arbitrary radial rescaling preserves Sobolev domains or finite-support density. For example, $a_0\equiv1$ and $a_1(x)=4^{|x|}$ satisfy $\mathfrak C_p(a_0)=\mathfrak C_p(a_1)=1$, whereas Proposition~\ref{prop:lp-equality} gives $W_{a_1}^{1,p}(T,m)\subsetneq W_{a_0}^{1,p}(T,m)=L^p(T,m)$. Outside the sufficient class $\{a:\mathfrak C_p(a)<\infty\}$, radial rescaling may also change finite-support density; see Remark~\ref{rem:rescaling-limit}. The same invariance identity need not hold for $\mathfrak C_p(1+a)$.
\end{remark}

In particular, if a positive, finite-valued radial coefficient $a(x)=A_{|x|}$ is unbounded, Proposition~\ref{prop:lp-equality} and Corollary~\ref{cor:radial} show that $W_a^{1,p}(T,m)$ is a proper subspace of $L^p(T,m)$ while $C_c(T)$ is still dense in $W_a^{1,p}(T,m)$.

For Examples~\ref{ex:branch-moment} and~\ref{ex:shift-enlarges}, consider a specific nonhomogeneous flow-measure tree. Fix an infinite ray $(r_n)_{n\ge0}$ starting at the root $\rho$, with $r_0=\rho$ and $\pi(r_n)=r_{n-1}$. Each $r_{n-1}$ has two children: the child $r_n$ receives the proportion
\[
 q_n=\begin{cases}1/2,&n\text{ odd},\\ 1/3,&n\text{ even},\end{cases}
\]
of the flow mass, while the other child receives $1-q_n$. Every vertex off the ray has three children, with flow proportions $(1/2,1/3,1/6)$. Normalize $m(\rho)=1$, and write $m_n=m(r_n)$. Then
\[
 m_{2j}=6^{-j},\qquad m_{2j+1}=\tfrac12\,6^{-j}\qquad(j\ge0).
\]
This defines a positive flow measure on a locally finite tree with nonuniform branching and flow splitting.

\begin{example}\label{ex:branch-moment}
On the nonhomogeneous flow-measure tree above, let $k(x)$ be the number of initial edges shared by the path from $\rho$ to $x$ and the distinguished ray. Let $0<\beta<1$, and let $(A_n)_{n\ge1}$ be any sequence of positive finite numbers. Define
\begin{equation*}
 a(x)=A_{|x|}2^{\beta k(x)},\qquad x\in\Ts.
\end{equation*}
Then
\[
 \frac{\max_{x\in S_n}a(x)}{\min_{x\in S_n}a(x)}
 =2^{\beta n}\longrightarrow\infty,
\]
while for every $1<p<\infty$,
\[
 \mathfrak C_p(a)\le\frac1{2-2^\beta}<\infty.
\]
Consequently, $C_c(T)$ is dense in $W_a^{1,p}(T,m)$.
\end{example}
\begin{proof}
At depth $n$, the largest and smallest values of $k(x)$ are $n$ and $0$, giving the stated coefficient ratio. To estimate the moment constant, first suppose that $z$ lies off the distinguished ray. Then $k$ is constant on $C(z,n)$, so $K_p(a;z,n)=1$.

If $z=r_N$, put $d=n-N\ge1$, $t=2^\beta\in(1,2)$, and $s_j=m_{N+j}/m_N$ for $1\le j\le d$. Since $q_n\le1/2$, we have $s_j\le2^{-j}$. Grouping descendants by the first edge leaving the ray (or by remaining on it through depth $n$) gives
\[
\begin{aligned}
 \frac{\av{a}{r_N,n}}{A_nt^N}
 &=1+(t-1)\sum_{j=1}^{d}t^{j-1}s_j\\
 &\le1+(t-1)\sum_{j=1}^{\infty}t^{j-1}2^{-j}
 =\frac1{2-t}.
\end{aligned}
\]
Since $k(x)\ge N$ on $C(r_N,n)$,
\[
 \av{a^{-1/(p-1)}}{r_N,n}^{p-1}\le(A_nt^N)^{-1}.
\]
Thus $K_p(a;r_N,n)\le(2-2^\beta)^{-1}$, uniformly in $N,n$. Together with the off-ray case, this proves the moment bound. Density now follows from Theorem~\ref{thm:density}.
\end{proof}

\begin{proposition}\label{prop:additive}
Let $a,b:\Ts\to(0,\infty)$ and suppose
\[
M=\sup_{x\in\Ts}|a(x)-b(x)|<\infty.
\]
Then, for every $f\in L^p(T,m)$,
\[
Q_{p,a}(f)<\infty
\iff
Q_{p,b}(f)<\infty.
\]
The corresponding Sobolev norms are equivalent, and $C_c(T)$ is dense in $W_a^{1,p}(T,m)$ if and only if it is dense in $W_b^{1,p}(T,m)$. In particular, every $f\in W_a^{1,p}(T,m)$ satisfies
\begin{equation}\label{eq:unit-shift}
\norm{f}_{W_a^{1,p}}^p
\le \norm{f}_{W_{1+a}^{1,p}}^p
\le (1+2^p)\norm{f}_{W_a^{1,p}}^p.
\end{equation}
\end{proposition}
\begin{proof}
By \eqref{eq:D-bound},
\[
Q_{p,a}(f)\le Q_{p,b}(f)+M Q_{p,1}(f)
\le Q_{p,b}(f)+2^pM\norm{f}_p^p.
\]
Interchanging $a$ and $b$ gives the reverse estimate. Thus the finite-energy conditions are equivalent. On the common domain,
\[
\begin{aligned}
\norm{f}_{W_a^{1,p}}^p
&\le (1+2^pM)\norm{f}_{W_b^{1,p}}^p,\\
\norm{f}_{W_b^{1,p}}^p
&\le (1+2^pM)\norm{f}_{W_a^{1,p}}^p.
\end{aligned}
\]
Hence the norms are equivalent and finite-support density is preserved. The additive error is controlled by the zero-order $L^p$ term. Taking $b=1+a$ and using $Q_{p,1+a}=Q_{p,a}+Q_{p,1}$ together with $Q_{p,1}(f)\le2^p\norm{f}_p^p$ yields \eqref{eq:unit-shift}.
\end{proof}

\begin{corollary}\label{cor:shifted-density}
Let $a:\Ts\to(0,\infty)$ be a positive, finite-valued coefficient. Then, for every $z\in T$ and $n>|z|$,
\begin{equation}\label{eq:shift-moment}
K_p(1+a;z,n)\le K_p(a;z,n),\qquad
\mathfrak C_p(1+a)\le\mathfrak C_p(a).
\end{equation}
If $\mathfrak C_p(1+a)<\infty$, then $C_c(T)$ is dense in $W_a^{1,p}(T,m)$.
\end{corollary}
\begin{proof}
Fix $z\in T$ and $n>|z|$. On $C(z,n)$, let $d\mu=dm/m(z)$ and let $\mathsf E$ denote the corresponding averaging operator. Set $K=K_p(a;z,n)\ge1$. By Jensen's inequality and Lemma~\ref{lem:weighted-average},
\[
\begin{aligned}
\norm{\mathsf Eh}_{L^p((1+a)d\mu)}^p
&=\norm{\mathsf Eh}_{L^p(d\mu)}^p
+\norm{\mathsf Eh}_{L^p(a\,d\mu)}^p\\
&\le \norm{h}_{L^p(d\mu)}^p
+K\norm{h}_{L^p(a\,d\mu)}^p\\
&\le K\norm{h}_{L^p((1+a)d\mu)}^p.
\end{aligned}
\]
The exact operator norm formula in Lemma~\ref{lem:weighted-average} therefore gives $K_p(1+a;z,n)\le K_p(a;z,n)$. Taking suprema yields \eqref{eq:shift-moment}. If $\mathfrak C_p(1+a)<\infty$, Theorem~\ref{thm:density} implies that $C_c(T)$ is dense in $W_{1+a}^{1,p}(T,m)$. The norm equivalence in Proposition~\ref{prop:additive} then gives density in $W_a^{1,p}(T,m)$.
\end{proof}

\begin{example}\label{ex:shift-enlarges}
On the nonhomogeneous flow-measure tree defined above, let $(r_n)_{n\ge0}$ be the distinguished ray and set
\[
 a(x)=\begin{cases}
 2^{-n^2},&x=r_n,\quad n\ge1,\\
 1,&x\in\Ts\setminus\{r_n:n\ge1\}.
 \end{cases}
\]
Then $\mathfrak C_p(a)=\infty$ and $\mathfrak C_p(1+a)\le2<\infty$ for every $1<p<\infty$. Consequently, $C_c(T)$ is dense in $W_a^{1,p}(T,m)$, and the shifted moment condition defines a strictly larger sufficient class than the original moment condition.
\end{example}
\begin{proof}
Put $\varepsilon_n=2^{-n^2}$. On $C(r_{n-1},n)$, the ray child and the other child have normalized flow masses $q_n$ and $1-q_n$, with coefficient values $\varepsilon_n$ and $1$, respectively. Thus
\[
\begin{aligned}
 K_p(a;r_{n-1},n)
 &=(q_n\varepsilon_n+1-q_n)
   (q_n\varepsilon_n^{-1/(p-1)}+1-q_n)^{p-1}\\
 &\ge(1-q_n)q_n^{p-1}\varepsilon_n^{-1}\\
 &\ge\frac{2^{n^2}}{2\cdot3^{p-1}}\longrightarrow\infty.
\end{aligned}
\]
On the other hand, $1\le1+a\le2$, whence $\mathfrak C_p(1+a)\le2$. The density conclusion follows from Corollary~\ref{cor:shifted-density}.
\end{proof}

\subsection{Limitations of the moment criterion and single-ray models}\label{subsec:moment-limitations}
The moment conditions above are sufficient rather than necessary for finite-support density. We now examine coefficients beyond these conditions using cutoff functions and limits along infinite rays.

To analyze density without the moment condition, we first record an alternative cutoff criterion.

\begin{lemma}\label{lem:small-energy-cutoff}
Suppose there exists a sequence $(\eta_N)_{N\ge0}\subset C_c(T)$ satisfying
\[
 0\le\eta_N\le1,\qquad \eta_N=1\ \text{on }B_N,
 \qquad Q_{p,a}(\eta_N)\longrightarrow0.
\]
Then $C_c(T)$ is dense in $W_a^{1,p}(T,m)$.
\end{lemma}
\begin{proof}
First let $f\in W_a^{1,p}(T,m)\cap\ell^\infty(T)$ and set $u_N=(1-\eta_N)f$. Since $u_N=0$ on $B_N$ and $0\le1-\eta_N\le1$,
\[
 \norm{u_N}_p^p\le\sum_{|x|>N}m(x)|f(x)|^p\longrightarrow0.
\]
Moreover,
\[
 Du_N(x)=(1-\eta_N(x))Df(x)-f(\pi(x))D\eta_N(x),
\]
and hence
\[
\begin{aligned}
 Q_{p,a}(u_N)
 &\le2^{p-1}\sum_{|x|>N}a(x)m(x)|Df(x)|^p\\
 &\quad+2^{p-1}\norm f_\infty^p Q_{p,a}(\eta_N)\longrightarrow0.
\end{aligned}
\]
Thus $\eta_Nf\to f$ in $W_a^{1,p}(T,m)$.

For general $f\in W_a^{1,p}(T,m)$, define $f_M=\max\{-M,\min\{f,M\}\}$. Both the truncation map and its residual are $1$-Lipschitz; the residual is bounded in absolute value by $|f|$ and converges pointwise to zero. Dominated convergence for the vertex and weighted-edge sums therefore gives $\norm{f-f_M}_{W_a^{1,p}}\to0$. For each fixed $M$, the bounded-function case gives $\norm{\eta_Nf_M-f_M}_{W_a^{1,p}}\to0$. Choose $M_j$ and then $N_j$ so that
\[
 \norm{f-f_{M_j}}_{W_a^{1,p}}<\frac1j,
 \qquad\norm{f_{M_j}-\eta_{N_j}f_{M_j}}_{W_a^{1,p}}<\frac1j.
\]
The finitely supported functions $g_j=\eta_{N_j}f_{M_j}$ satisfy $\norm{f-g_j}_{W_a^{1,p}}<2/j\to0$, proving density.
\end{proof}

\begin{proposition}\label{prop:ray}
Let $(r_n)_{n\ge0}$ be an infinite ray with $r_0=\rho$ and $\pi(r_n)=r_{n-1}$ for $n\ge1$. Suppose
\begin{align}
 R&=\sum_{n=1}^{\infty}[a(r_n)m(r_n)]^{-1/(p-1)}<\infty,\notag\\
 M&=\sum_{n=0}^{\infty}m(r_n)<\infty,\label{eq:ray-mass}\\
 S&=\sum_{n=0}^{\infty}\sum_{\substack{y\in s(r_n)\\y\ne r_{n+1}}}a(y)m(y)<\infty.\label{eq:ray-leak}
\end{align}
Then $W_{a,0}^{1,p}(T,m)\subsetneq W_a^{1,p}(T,m)$.
\end{proposition}
\begin{proof}
For $u\in W_a^{1,p}(T,m)$, weighted H\"older's inequality along the ray gives
\begin{equation*}
 \sum_{n=1}^{\infty}|Du(r_n)|
 \le\left(\sum_{n=1}^{\infty}a(r_n)m(r_n)|Du(r_n)|^p\right)^{1/p}R^{1/p'}
 \le Q_{p,a}(u)^{1/p}R^{1/p'}.
\end{equation*}
Thus the ray limit
\[
 \ell(u)=\lim_{n\to\infty}u(r_n)
 =u(\rho)+\sum_{n=1}^{\infty}Du(r_n)
\]
exists and defines a linear functional. Since $m(\rho)=1$, H\"older's inequality yields
\begin{equation*}
\begin{aligned}
 |\ell(u)|
 &\le|u(\rho)|+R^{1/p'}Q_{p,a}(u)^{1/p}\\
 &\le(1+R)^{1/p'}\norm u_{W_a^{1,p}}.
\end{aligned}
\end{equation*}
Hence $\ell$ is continuous. Since every finitely supported function has ray limit zero, $\ell$ vanishes on $W_{a,0}^{1,p}(T,m)$.

Now let $f=\ind_{\{r_0,r_1,\ldots\}}$. By \eqref{eq:ray-mass} and \eqref{eq:ray-leak},
\[
 \norm f_p^p=M,\qquad Q_{p,a}(f)=S.
\]
Thus $f\in W_a^{1,p}(T,m)$ and $\ell(f)=1$. Consequently, for every $v\in W_{a,0}^{1,p}(T,m)$,
\[
 1=|\ell(f-v)|\le(1+R)^{1/p'}\norm{f-v}_{W_a^{1,p}}.
\]
Taking the infimum yields
\begin{equation*}
 \dist_{W_a^{1,p}}\bigl(f,W_{a,0}^{1,p}(T,m)\bigr)
 \ge(1+R)^{-1/p'}>0,
\end{equation*}
which proves the strict inclusion.
\end{proof}

We return to the general assumptions on $(T,m)$ from Subsection~\ref{sec:prelim}. Fix an infinite ray $(r_n)_{n\ge0}$, with $r_0=\rho$ and $\pi(r_n)=r_{n-1}$. Write $m_n=m(r_n)$, and suppose there exist constants $0<q_-\le q_+<1$ such that
\[
 q_-\le q_n:=\frac{m_n}{m_{n-1}}\le q_+,\qquad n\ge1.
\]
For $\gamma>0$, define the single-ray coefficient
\[
 a_\gamma(x)=\begin{cases}
 m_n^{-\gamma},&x=r_n,\quad n\ge1,\\
 1,&x\in\Ts\setminus\{r_n:n\ge1\}.
 \end{cases}
\]
The concrete nonhomogeneous tree described in Subsection~\ref{subsec:moment-applications} satisfies these assumptions with $q_-=1/3$ and $q_+=1/2$.

\begin{theorem}\label{thm:ray-family}
For every $1<p<\infty$ and $\gamma>0$,
\begin{equation*}
 \mathfrak C_p(a_\gamma)=\mathfrak C_p(1+a_\gamma)=\infty,
 \qquad W_{a_\gamma}^{1,p}(T,m)\subsetneq L^p(T,m).
\end{equation*}
Moreover, $C_c(T)$ is dense in $W_{a_\gamma}^{1,p}(T,m)$ if and only if
\begin{equation*}
 0<\gamma\le1.
\end{equation*}
\end{theorem}
\begin{proof}
Since $m_0=1$ and $0<q_n\le q_+<1$,
\[
 0<m_n=\prod_{j=1}^{n}q_j\le q_+^n\longrightarrow0.
\]
We first show that both moment constants are infinite. Put $t_n=m_n^{-\gamma}$ and $c=q_-(1-q_+)^{p-1}>0$. On $C(r_{n-1},n)$, the ray child has normalized flow mass $q_n$ and coefficient $t_n$; the other children have total normalized mass $1-q_n$ and coefficient $1$. Thus
\[
\begin{aligned}
 K_p(a_\gamma;r_{n-1},n)
 &=(q_nt_n+1-q_n)
 (q_nt_n^{-1/(p-1)}+1-q_n)^{p-1}\\
 &\ge q_n(1-q_n)^{p-1}t_n
 \ge ct_n\longrightarrow\infty.
\end{aligned}
\]
For $1+a_\gamma$, the corresponding coefficient values are $1+t_n$ and $2$, hence
\[
\begin{aligned}
 K_p(1+a_\gamma;r_{n-1},n)
 &\ge q_nt_n\bigl((1-q_n)2^{-1/(p-1)}\bigr)^{p-1}\\
 &\ge\tfrac c2 t_n\longrightarrow\infty.
\end{aligned}
\]
Thus both moment constants are infinite. Since $a_\gamma(r_n)=t_n\to\infty$, Proposition~\ref{prop:lp-equality} gives $W_{a_\gamma}^{1,p}(T,m)\subsetneq L^p(T,m)$.

We now prove the density criterion. If $0<\gamma\le1$, the flow identity gives, for every $n\ge1$,
\begin{equation*}
 \sum_{x\in S_n}a_\gamma(x)m(x)=m_n^{1-\gamma}+1-m_n\le2.
\end{equation*}
For $N\ge1$, choose $\eta_N(x)=\theta_N(|x|)$, where
\[
 \theta_N(n)=\begin{cases}
 1,&n\le N,\\
 (2N-n)/N,&N<n<2N,\\
 0,&n\ge2N.
 \end{cases}
\]
Local finiteness gives $\eta_N\in C_c(T)$, and $0\le\eta_N\le1$, $\eta_N=1$ on $B_N$. Its nonzero differences occur at depths $N+1,\ldots,2N$, each with absolute value $1/N$. Hence
\[
 Q_{p,a_\gamma}(\eta_N)
 =N^{-p}\sum_{n=N+1}^{2N}\sum_{x\in S_n}a_\gamma(x)m(x)
 \le2N^{1-p}\longrightarrow0.
\]
By Lemma~\ref{lem:small-energy-cutoff}, finite-support density holds, including $\gamma=1$.

Conversely, if $\gamma>1$, the quantities in Proposition~\ref{prop:ray} satisfy
\[
 R=\sum_{n=1}^{\infty}m_n^{(\gamma-1)/(p-1)}
 \le\sum_{n=1}^{\infty}q_+^{\,n(\gamma-1)/(p-1)}<\infty,
\]
\[
 M=\sum_{n=0}^{\infty}m_n\le\frac1{1-q_+}<\infty.
\]
Since $a_\gamma=1$ on every edge leaving the distinguished ray, flow conservation gives
\[
 S=\sum_{n=0}^{\infty}(m_n-m_{n+1})=1.
\]
Proposition~\ref{prop:ray} then proves that $C_c(T)$ is not dense in $W_{a_\gamma}^{1,p}(T,m)$. The stated equivalence follows.
\end{proof}

\begin{example}\label{ex:same-max}
Retain the tree, flow measure, and distinguished ray of Theorem~\ref{thm:ray-family}, and write $m_n=m(r_n)$. Set
\[
 a_{\rm lev}(x)=m_{|x|}^{-2},\qquad a_{\rm ray}=a_2,\qquad x\in\Ts,
\]
where $a_2$ is the single-ray coefficient with $\gamma=2$. For every $n\ge1$,
\begin{equation*}
\begin{aligned}
 \max_{x\in S_n}a_{\rm lev}(x)
 &=\max_{x\in S_n}a_{\rm ray}(x)=m_n^{-2},\\
 a_{\rm lev}(r_n)&=a_{\rm ray}(r_n)=m_n^{-2}.
\end{aligned}
\end{equation*}
However, the off-ray edge weights differ. For the ray indicator $f=\ind_{\{r_0,r_1,\ldots\}}$, we have
\[
 \norm f_p^p=\sum_{n=0}^{\infty}m_n<\infty,
\]
and flow conservation yields
\[
\begin{aligned}
 Q_{p,a_{\rm ray}}(f)
 &=\sum_{n=1}^{\infty}(m_{n-1}-m_n)=1,\\
 Q_{p,a_{\rm lev}}(f)
 &=\sum_{n=1}^{\infty}m_n^{-2}(m_{n-1}-m_n)=\infty.
\end{aligned}
\]
Indeed, each summand in the second series satisfies
\[
 m_n^{-2}(m_{n-1}-m_n)
 =\frac{1-q_n}{q_nm_n}
 \ge\frac{1-q_+}{q_+m_n}\longrightarrow\infty.
\]
By Corollary~\ref{cor:radial} and Theorem~\ref{thm:ray-family}, $C_c(T)$ is dense in $W_{a_{\rm lev}}^{1,p}(T,m)$ but not in $W_{a_{\rm ray}}^{1,p}(T,m)$. Thus neither levelwise maxima nor coefficient values along a distinguished ray alone determine finite-support density; the distribution of coefficients across different branches must also be taken into account.
\end{example}

\begin{remark}\label{rem:rescaling-limit}
With the notation of Example~\ref{ex:same-max}, let
\[
 b(x)=m_{|x|}^{2}a_{\rm ray}(x),\qquad x\in\Ts.
\]
Since $0<b\le1$, Proposition~\ref{prop:lp-equality} implies that $C_c(T)$ is dense in $W_b^{1,p}(T,m)$, whereas Theorem~\ref{thm:ray-family} shows that it is not dense in $W_{a_{\rm ray}}^{1,p}(T,m)$. Nevertheless, Proposition~\ref{prop:radial-invariance} gives
\[
 \mathfrak C_p(b)=\mathfrak C_p(a_{\rm ray})=\infty.
\]
The effective edge weights along the distinguished ray are
\[
 b(r_n)m(r_n)=m_n,\qquad
 a_{\rm ray}(r_n)m(r_n)=m_n^{-1}.
\]
Thus radial rescaling preserves the weighted moment constant but may change the energy constraints along a ray and hence finite-support density.
\end{remark}

\section{Applications}\label{sec:applications}

\subsection{A diffusion model on a branching medium}\label{subsec:diffusion}
Consider a stationary diffusion model on the tree $T$. Let $u:T\to(-\infty,\infty)$ denote the state at each vertex, and set
\[
 c_x=a(x)m(x),\qquad x\in\Ts,
\]
as the effective response coefficient of the corresponding edge. For $1<p<\infty$, prescribe the edge flux
\[
 j_x=c_x|Du(x)|^{p-2}Du(x),
\]
with the positive direction taken from $x$ to its parent $\pi(x)$, and choose the local vertex response at $z$ to be
\[
 m(z)|u(z)|^{p-2}u(z).
\]
Here $|s|^{p-2}s$ is defined to be zero at $s=0$. Taking $m(z)h(z)$ as the external source and setting $j_\rho=0$, the stationary balance equation is
\begin{equation*}
 m(z)|u(z)|^{p-2}u(z)+j_z-\sum_{x\in s(z)}j_x=m(z)h(z).
\end{equation*}
When $p=2$, this reduces to a linear diffusion model. In this case, $u$ may represent the temperature deviation from a fixed reference temperature, $c_x$ may be interpreted as a normalized effective edge conductance, and the vertex response describes linear exchange with the reference environment. For $p\ne2$, the edge flux follows a nonlinear power law rather than the classical Fourier conduction law.

The flow measure $m$ serves as a reference measure, not as a physical heat flux. These edge and vertex responses correspond, respectively, to the weighted difference term and the zero-order term in the Sobolev norm defined in Subsection~\ref{sec:prelim}. The choice of the corresponding weights is a modeling assumption, not a consequence of the flow identity. For related variational theory of discrete $p$-Laplacians, see \cite{Mugnolo2013}. The local operator and variational formulation for the present model are given in Subsection~\ref{subsec:weak-operator}.

\subsection{Variational formulation and compactly supported testing}\label{subsec:weak-operator}
By Lemma~\ref{lem:D-bound}, for $p=2$, the difference operator
\[
 D:L^2(T,m)\longrightarrow L^2(\Ts,m)
\]
is bounded. Its Hilbert-space adjoint $D^*:L^2(\Ts,m)\to L^2(T,m)$ is given by
\begin{equation}\label{eq:adjoint}
 D^*g(z)=
 \begin{cases}
 g(z)-\displaystyle\frac1{m(z)}\sum_{x\in s(z)}m(x)g(x), &z\ne\rho,\\[7pt]
 -\displaystyle\frac1{m(\rho)}\sum_{x\in s(\rho)}m(x)g(x), &z=\rho.
 \end{cases}
\end{equation}
Indeed, summation by parts gives
\[
 \sum_{x\in\Ts}Df(x)g(x)m(x)=\sum_{z\in T}f(z)D^*g(z)m(z)
\]
for finitely supported $f$ and $g$. By density and boundedness, this identity extends to the corresponding $L^2$ spaces.

When $a\equiv1$, the local expression of $\frac12D^*D$ at nonroot vertices agrees with the normalized flow Laplacian in \cite{LMSTV2023,MSTV2025}. These references consider trees rooted at infinity, whereas the distinguished root $\rho$ in the present setting has no parent edge. Consequently, the expression at $\rho$ is different.

With the stationary flux prescribed in Subsection~\ref{subsec:diffusion}, define the local weighted discrete $p$-Laplacian, for any real-valued function $u:T\to(-\infty,\infty)$ and $z\in\Ts$, by
\begin{equation*}
\begin{split}
 \mathscr A_{p,a}u(z)
 &=a(z)|u(z)-u(\pi(z))|^{p-2}\bigl(u(z)-u(\pi(z))\bigr)\\
 &\quad+\frac1{m(z)}\sum_{x\in s(z)}a(x)m(x)
 |u(z)-u(x)|^{p-2}\bigl(u(z)-u(x)\bigr).
\end{split}
\end{equation*}
At the root $\rho$, only the second term is present. Since $T$ is locally finite and the coefficients are positive and finite-valued, this expression is well-defined at every vertex.

Using the local divergence formula \eqref{eq:adjoint}, we may formally write
\begin{equation}\label{eq:formal-p}
 \mathscr A_{p,a}u=D^*\bigl(a|Du|^{p-2}Du\bigr).
\end{equation}
Here $D^*$ is understood pointwise through the divergence formula, and \eqref{eq:formal-p} does not in general represent a composition of operators on $L^2$. For $p=2$, it is a genuine composition with the Hilbert-space adjoint whenever $aDu\in L^2(\Ts,m)$. However, membership in $W_a^{1,2}(T,m)$ guarantees only $a^{1/2}Du\in L^2(\Ts,m)$ and does not imply this additional condition.

For $u,\varphi\in\Wa$, H\"older's inequality gives
\begin{equation}\label{eq:variation-bound}
 \left|\sum_{x\in\Ts}a(x)m(x)|Du(x)|^{p-2}Du(x)D\varphi(x)\right|
 \le\Qa(u)^{1/p'}\Qa(\varphi)^{1/p}.
\end{equation}
Thus the sum converges absolutely and gives the first variation of $Q_{p,a}/p$:
\[
 \left.\frac{d}{dt}\frac{Q_{p,a}(u+t\varphi)}p\right|_{t=0}
 =\sum_{x\in\Ts}a(x)m(x)|Du(x)|^{p-2}Du(x)D\varphi(x).
\]
For $\varphi\in\Cc$, the difference $D\varphi$ also has finite support, and summation by parts yields
\begin{equation}\label{eq:local-green}
 \sum_{z\in T}\mathscr A_{p,a}u(z)\varphi(z)m(z)
 =\sum_{x\in\Ts}a(x)m(x)|Du(x)|^{p-2}Du(x)D\varphi(x),
 \qquad \varphi\in\Cc.
\end{equation}
No boundary condition at infinity is required for this identity.

For $h\in L^{p'}(T,m)$, define on $\Wa$ the functional
\begin{equation*}
 \mathcal J_h(u)=\frac1p\left(\norm u_p^p+\Qa(u)\right)
 -\sum_{z\in T}m(z)h(z)u(z).
\end{equation*}
We consider its minimization over $\Wa$ and $\Waz$. The corresponding local stationary equation is
\begin{equation*}
 |u|^{p-2}u+\mathscr A_{p,a}u=h.
\end{equation*}
For $u,\varphi\in\Wa$, the first variation of $\mathcal J_h$ is
\begin{equation*}
\begin{split}
 \mathcal J_h'(u)[\varphi]
 &=\sum_{z\in T}m(z)|u(z)|^{p-2}u(z)\varphi(z)\\
 &\quad+\sum_{x\in\Ts}a(x)m(x)|Du(x)|^{p-2}Du(x)D\varphi(x)\\
 &\quad-\sum_{z\in T}m(z)h(z)\varphi(z).
\end{split}
\end{equation*}
By \eqref{eq:variation-bound} and H\"older's inequality,
\[
\begin{split}
 |\mathcal J_h'(u)[\varphi]|
 &\le\bigl(\norm u_p^{p-1}+\norm h_{p'}\bigr)\norm\varphi_p
 +\Qa(u)^{1/p'}\Qa(\varphi)^{1/p}.
\end{split}
\]
Hence, for fixed $u$, the functional $\mathcal J_h'(u)$ is continuous and linear on $\Wa$. If it vanishes on all finitely supported test functions, continuity implies that it also vanishes on $\Waz$. When $\Cc$ is dense in $\Wa$, the identity extends to every Sobolev test function.

Moreover,
\[
 \mathcal J_h(u)\ge\frac1p\norm u_{\Wa}^p-\norm h_{p'}\norm u_{\Wa}.
\]
Since $1<p<\infty$, the functional is coercive, strictly convex, and weakly lower semicontinuous on both energy domains. Together with Theorem~\ref{thm:structure}, these properties ensure a unique minimizer on each domain, whose first variation vanishes on the corresponding test space.

For $p=2$, define the quadratic form
\[
 q_a(u,v)=\sum_{x\in\Ts}a(x)m(x)Du(x)Dv(x)
\]
with form domain $W_a^{1,2}(T,m)$. Since $\Cc$ is dense in $L^2(T,m)$ and the form norm satisfies
\[
 \norm u_{q_a}^2=\norm u_2^2+q_a(u,u)=\norm u_{W_a^{1,2}}^2,
\]
Theorem~\ref{thm:structure} shows that $q_a$ is a densely defined, nonnegative, closed quadratic form. Its associated nonnegative self-adjoint operator is characterized as follows: $u\in W_a^{1,2}(T,m)$ belongs to the operator domain if and only if there exists $g\in L^2(T,m)$ such that
\[
 q_a(u,v)=\langle g,v\rangle_{L^2(T,m)}\qquad\text{for every }v\in W_a^{1,2}(T,m).
\]
The operator then acts as $g=\mathscr A_{2,a}u$ pointwise. Using $W_{a,0}^{1,2}(T,m)$ as the form domain instead gives the self-adjoint realization associated with the finite-support closure. The two realizations have the same local expression but different form domains and operator realizations when finitely supported functions are not dense; see \cite{HKLW2012} for the related theory of quadratic forms on graphs. Subsection~\ref{subsec:quadratic} discusses the consequences for variational solutions and finite-domain approximation.

\subsection{The quadratic case: finite-domain and global solutions}\label{subsec:quadratic}
Let
\[
 H=W_a^{1,2}(T,m),\qquad H_0=W_{a,0}^{1,2}(T,m),
\]
and define
\begin{equation*}
 B_a(u,v)=\sum_{z\in T}u(z)v(z)m(z)+\sum_{x\in\Ts}a(x)m(x)Du(x)Dv(x),
 \quad B_a(u,u)=\norm u_H^2.
\end{equation*}
Since $H_0$ is a closed subspace of $H$, let $\Pi_0:H\to H_0$ be the $B_a$-orthogonal projection, and set
\[
 \mathcal K_a=H_0^{\perp_{B_a}}
 =\{w\in H:B_a(w,v)=0\text{ for every }v\in H_0\}.
\]
Then $H=H_0\oplus_{B_a}\mathcal K_a$.

\begin{proposition}\label{prop:quadratic}
For each $h\in L^2(T,m)$, there exist unique $u^{\max}\in H$ and $u^0\in H_0$ satisfying
\begin{equation}\label{eq:two-weak-problems}
\begin{aligned}
 B_a(u^{\max},v)&=\langle h,v\rangle_{L^2(T,m)} &&(v\in H),\\
 B_a(u^0,v)&=\langle h,v\rangle_{L^2(T,m)} &&(v\in H_0).
\end{aligned}
\end{equation}
Moreover, $\norm{u^{\max}}_H,\ \norm{u^0}_H\le\norm h_2$ and $u^0=\Pi_0u^{\max}$. All $u\in H$ satisfying
\begin{equation}\label{eq:local-weak}
 B_a(u,\varphi)=\langle h,\varphi\rangle_{L^2(T,m)},\qquad \varphi\in\Cc
\end{equation}
form the affine space $u^0+\mathcal K_a$. Consequently, compactly supported testing determines a unique finite-energy solution in $H$ if and only if $\Cc$ is dense in $H$. If density fails, the two solutions $u^{\max}$ and $u^0$ differ for some $h\in L^2(T,m)$.
\end{proposition}
\begin{proof}
Since $|\langle h,v\rangle_{L^2(T,m)}|\le\norm h_2\norm v_H$, the Riesz representation theorem yields the unique solutions of \eqref{eq:two-weak-problems} and the stated norm bounds. Subtracting the two variational identities on $H_0$ gives $u^{\max}-u^0\in\mathcal K_a$, and hence $u^0=\Pi_0u^{\max}$.

On the other hand, \eqref{eq:local-weak} is equivalent to $B_a(u-u^0,\varphi)=0$ for every $\varphi\in\Cc$. By the density of $\Cc$ in $H_0$ and continuity of $B_a$, this holds precisely when $u-u^0\in\mathcal K_a$. Thus the solution set is $u^0+\mathcal K_a$, and uniqueness in $H$ holds exactly when $\mathcal K_a=\{0\}$, or equivalently, when $\Cc$ is dense in $H$.

Finally, if $\mathcal K_a\ne\{0\}$, choose $0\ne w\in\mathcal K_a$ and set $h=w$. Then $B_a(u^{\max},w)=\norm w_2^2>0$, whereas $B_a(u^0,w)=0$. Therefore $u^{\max}\ne u^0$.
\end{proof}

Testing against vertex indicators and using \eqref{eq:local-green} shows that \eqref{eq:local-weak} is equivalent to the pointwise equation $u+\mathscr A_{2,a}u=h$. In particular, $\mathcal K_a$ consists of finite-energy solutions of $w+\mathscr A_{2,a}w=0$. The existence of a nonzero such solution does not contradict the positive definiteness of $B_a$, since compactly supported testing does not allow $w\notin H_0$ itself as a test function.

For $N\ge0$, let
\[
 V_N=\{v\in H:\supp v\subseteq B_N\}.
\]
By local finiteness, $V_N$ is finite-dimensional, and $V_N\subseteq H_0$ with $\bigcup_{N\ge0}V_N=\Cc$.

\noindent\textbf{Proof of Theorem~\ref{thm:finite-domain-convergence}.} We now prove Theorem~\ref{thm:finite-domain-convergence}.
\begin{proof}
Fix $h\in L^2(T,m)$. The finite-domain solution $u_N\in V_N$ is uniquely determined by
\begin{equation}\label{eq:finite-weak}
 B_a(u_N,v)=\langle h,v\rangle_{L^2(T,m)},\qquad v\in V_N.
\end{equation}
Subtracting \eqref{eq:finite-weak} from the variational identity for $u^0$ in Proposition~\ref{prop:quadratic} yields $B_a(u^0-u_N,v)=0$ for all $v\in V_N$. Thus $u_N$ is the $B_a$-orthogonal projection of $u^0$ onto $V_N$. Since the spaces $V_N$ are increasing and their union is dense in $H_0$,
\begin{equation}\label{eq:finite-best}
 \norm{u^0-u_N}_H=\inf_{v\in V_N}\norm{u^0-v}_H\longrightarrow0.
\end{equation}
Moreover, $u^{\max}-u^0\in\mathcal K_a$ and $u^0-u_N\in H_0$. These differences are orthogonal, giving
\begin{equation}\label{eq:finite-gap}
 \norm{u^{\max}-u_N}_H^2
 =\norm{u^{\max}-u^0}_H^2+\norm{u^0-u_N}_H^2.
\end{equation}
If $\Cc$ is dense in $H$, then $u^{\max}=u^0$, and \eqref{eq:finite-best} implies $u_N\to u^{\max}$. Conversely, if density fails, Proposition~\ref{prop:quadratic} provides some $h$ for which $u^{\max}\ne u^0$. Equation~\eqref{eq:finite-gap} then excludes convergence to $u^{\max}$ for that datum. This proves the theorem.
\end{proof}

For a fixed datum $h$, \eqref{eq:finite-best}--\eqref{eq:finite-gap} show that $u_N\to u^{\max}$ if and only if $u^{\max}=u^0$. The equivalence in Theorem~\ref{thm:finite-domain-convergence} concerns convergence for every $h\in L^2(T,m)$.

Finally, consider the finite linear system associated with the variational problem. Write $y\sim z$ when $y$ and $z$ are adjacent. Retaining the effective edge weights $c_x=a(x)m(x)$ from Subsection~\ref{subsec:diffusion}, define the symmetric coefficients $c_{\pi(x),x}=c_{x,\pi(x)}=c_x$ on each edge. Extending $u_N$ by zero outside $B_N$ and testing \eqref{eq:finite-weak} against vertex indicators gives
\begin{equation*}
 \left(m(z)+\sum_{y\sim z}c_{zy}\right)u_N(z)
 -\sum_{\substack{y\sim z\\y\in B_N}}c_{zy}u_N(y)
 =m(z)h(z),\qquad z\in B_N.
\end{equation*}
The coefficient matrix is symmetric positive definite, since its quadratic form is the restriction of $B_a$ to $V_N$. The diagonal terms retain edges crossing the boundary of $B_N$. Thus the finite-domain problem imposes a zero exterior state rather than an insulating boundary condition.

If $\mathfrak C_2(a)<\infty$, or more generally $\mathfrak C_2(1+a)<\infty$, Theorem~\ref{thm:density} and Corollary~\ref{cor:shifted-density} imply that $\Cc$ is dense in $H$. Consequently, the finite-domain solutions converge to $u^{\max}$ for every $h\in L^2(T,m)$.

For the single-ray coefficients $a_\gamma$ in Theorem~\ref{thm:ray-family}, the finite-domain solutions converge to $u^{\max}$ for every $h\in L^2(T,m)$ if and only if $0<\gamma\le1$, even though $\mathfrak C_2(a_\gamma)=\mathfrak C_2(1+a_\gamma)=\infty$.

\bigskip

\noindent{\bf\Large Declarations}

\noindent{\bf Acknowledgements } 
The authors would like to thank the editors and reviewers for careful reading and valuable comments, which lead to the improvement of this paper.

\medskip
\noindent{\bf Mathematics Subject Classification(2020)}
Primary 46E35; Secondary 05C05; 35R02; 35A15.

\medskip
\noindent{\bf Data, materials, and code availability}
Our manuscript has no associated data, materials, and code.

\medskip 
\noindent{\bf Funding information} The research was supported by the National Natural Science Foundation of China (No. 12461021).

\medskip
\noindent{\bf Contributions} All authors participated in the conception and design of the study and reviewed the manuscript.

\medskip
\noindent{\bf Conflict of interest} The authors declare that they have no conflict of interest.

\medskip
\noindent{\bf Author contributions}
All authors participated in the conception and design of the study and reviewed the manuscript.

\medskip
\noindent{\bf Consent for publication}
The authors agree to the publication.

\medskip
\noindent{\bf Ethical statement}
{\bf (1) Conflict of Interest:} The authors declare that they have no known competing financial interests or personal relationships that could have appeared to influence the work reported in this paper. {\bf (2) Data Availability:} No new data were generated or analyzed in support of this research. {\bf (3) Authorship and Originality:} This manuscript is the authors' original work and has not been published or submitted elsewhere. All authors have read and approved the final version and agree to be accountable for all aspects of the work.

\end{document}